\documentclass[aip,jcp]{revtex4-1}

\usepackage{amsmath,amssymb}
\usepackage{dcolumn}
\usepackage{bm}

\usepackage{amsthm}
\usepackage{mathrsfs}
\usepackage{graphicx}
\usepackage{hyperref}
\usepackage{xcolor}
\usepackage{fourier}
\usepackage{geometry}
\usepackage{extarrows}

\newcommand{\bl}{\big\langle}
\newcommand{\br}{\big\rangle}

\def\bl{\Big\langle}
\def\br{\Big\rangle}
\def\blue{\color{black}}

\begin{document}
\title{From Generalized Langevin Equations to Brownian Dynamics and Embedded Brownian Dynamics}

\author{Lina Ma}
\email{linama@psu.edu}
\affiliation{Department of Mathematics, the Pennsylvania State University, University Park, PA 16802-6400, USA.}%
\author{ Xiantao Li}
\email{xli@math.psu.edu}
\affiliation{Department of Mathematics, the Pennsylvania State University, University Park, PA 16802-6400, USA.}%
\author{Chun Liu}
\email{liu@math.psu.edu}
\affiliation{Department of Mathematics, the Pennsylvania State University, University Park, PA 16802-6400, USA.}%
\begin{abstract}
We present the reduction of generalized Langevin equations to a coordinate-only stochastic model, which in its exact form, involves a forcing term with memory and a general Gaussian noise. It will be shown that a similar fluctuation-dissipation theorem still holds at this level.  We study the approximation by the typical Brownian dynamics as a first approximation. Our numerical test indicates how the intrinsic frequency of the kernel function influences the accuracy of this approximation. In the case when such an approximate is inadequate,  further approximations can be derived by embedding the nonlocal model into an extended dynamics without memory.  By imposing noises in the auxiliary variables, we show how the second fluctuation-dissipation theorem is still exactly satisfied. 
\end{abstract}
\maketitle

\section{Introduction}

The Langevin dynamics (LD) model plays a crucial role in the stochastic modeling of  bio-molecules \cite{Schlick2002}. In general, the LD model can be expressed as follows,
\begin{equation}\label{eq: ld}
 m \ddot{x} = f(x) -\gamma \dot{x} + \xi(t). 
\end{equation}
Here $m$ denotes the mass, $f(x)$ is the inter-molecular force, and $\gamma$ is the friction coefficient.   In addition, $\xi(t)$ is an external white noise, which satisfies the fluctuation-dissipation theorem (FDT) \cite{Kubo66},
\begin{equation}\label{eq: fdt}
  \bl \xi(t) \xi(t')^T\br= 2k_B T \gamma \delta(t-t').
\end{equation}
This condition ensures that the dynamical system settles to the correct equilibrium \cite{Kubo66}. 

{\blue In the Langevin dynamics model \eqref{eq: ld}, the damping coefficient and the random noise are modeling the influence of the solvent particles.}  One interesting regime is where $\gamma \gg 1.$  The asymptotic analysis has been
the subject of theoretical interest.  The analysis can be done either through the Fokker-Planck equation \cite{titulaer1980corrections, titulaer1978systematic} using the Chapman-Enskog expansion, or directly based on
the stochastic differential equation \cite{san1980colored}. In {\blue Ref. \cite{san1980colored}},  the LD model was written in the first order form, and by solving the second equation, the velocity variables can be eliminated.  These analysis confirms that, when the damping is large, the inertial term can be neglected, leading to,
\begin{equation}\label{eq: bd}
  \gamma \dot{x}= f(x) + \xi(t). 
\end{equation}
This model has been referred by some authors as the Brownian dynamics (BD) model \cite{huber2010browndye,ricci2003algorithms,van1982algorithms}. Notice that in this reduced model, $\xi(t)$ still obeys
the FDT \eqref{eq: fdt}.

The derivation of the Brownian dynamics is of obvious theoretical interest. In practice, the BD (\ref{eq: bd}) neglects the transient time scale, and as a result, the dynamics occurs on a much larger time scale. Therefore, from a computational viewpoint, the time step can be much larger than that of the original Langevin dynamics (\ref{eq: ld}). This model is of considerable advantage in bio-molecular simulations, in which the time scale of interest is often out of reach when the full model \eqref{eq: ld} is used. 

There have been important theoretical works toward the asymptotic analysis in
the limit of large damping coefficients. The question of high friction limit also arises in
the study of plasma gases, especially those involving nonlocal Coulomb
interactions,  with kinetic descriptions \cite{el2010diffusion,goudon2005hydrodynamic,goudon2005multidimensional,  poupaud2000parabolic, nieto2001high}.
 In our recent work \cite{wu2015diffusion},
we studied a
Vlasov-Poisson-Fokker-Planck (VPFP) system in a bounded domain with reflection boundary
conditions for charge distributions, showed the convergence to the
macroscopic Poisson-Nerest-Planck (PNP) system. The purpose there was trying to
justify the PNP system as a macroscopic model for the transport of multi-species ions
in dilute solutions.

{\blue The BD model \eqref{eq: bd} can be viewed as a first order approximation to the full Langevin dynamics model \eqref{eq: bd}. } Meanwhile, there have been numerous attempts to derive high order approximations to improve the modeling accuracy without having to resolve the small time scales. A remarkable success is the inertial Brownian dynamics model (IBD), developed by Beard and Schlick \cite{beard2000inertial,beard2000inertialB}. Using further asymptotic expansions, the IBD incorporates correction terms into the BD model and takes into account the transient part of the dynamics. 

Motivated by these effort,  we study the BD approximation of the generalized Langevin equation (GLE). In contrast to the Langevin dynamics \eqref{eq: ld}, the GLE includes a history-dependent friction term and a correlated Gaussian noise. The memory term often comes from a coarse-graining step, e.g.,  \cite{AdDo74,tully1980dynamics,guardia1985generalized,ciccotti1981derivation,LiXian2014,ChSt05,curtarolo_dynamics_2002,espanol2004statistical,EvMo08,IzVo06,Li2009c,Li14,LiE07,li2015incorporation,MaLiLiu16,oliva2000generalized,stepanova_dynamics_2007,Zwanzig73,lei2016generalized}, which has also been an recent emerging area of interest in molecular modeling.  Among the many approximation schemes \cite{berkowitz1983generalized,berkowitz1981memory,LiXian2014,Darve_PNAS_2009,hijon2010mori,kauzlaric2011bottom,MaLiLiu16,oliva2000generalized} for the GLE
 model, the BD model (\ref{eq: bd}) is clearly the simplest. To understand how the BD approximation can naturally come about, we first present a derivation which eliminates momentum variables from the GLE, so that the resulting model only involves the coordinates of the molecular variables. It will be proved that the reduced model still exhibits history-dependence, and the second fluctuation-dissipation theorem, an important property of the GLE, is inherited by the reduced model, just like the BD approximation of the full Langevin dynamics. 
 
 At this level, the BD approximation arises as the first approximation, in which the memory kernel is approximated by a Dirac delta  function.
 To understand the accuracy and possibly, the limitation of this approximation, we consider the GLE derived by Adelman and Doll \cite{AdDo74,AdDo76}, and conducted several numerical tests. Our results indicate that the accuracy hinges on the frequency associated with the memory kernel: When the frequency is large, compared to that of the mean force, the approximation by BD is quite promising. Otherwise, the accuracy is quite limited and high order approximations are needed. To improve the modeling accuracy, we derive two other approximations where the memory kernel is approximated by rational functions in the Laplace domain.
The main motivation has been two-fold: On one hand, in the time domain, the memory can be eliminated by introducing auxiliary variables. As a result, there is no need to keep the history of the solution and evaluate the integral associated with the memory term at every step. This significantly suppresses  the computational cost. Meanwhile, it will be shown that with proper choices of the stochastic forces, applied to the auxiliary variables, the noise in the GLE is approximated by a colored noise, in such a way that the second FDT is exactly satisfied. 

\section{Mathematical Derivation}\label{eq: sec-deri}

We start with the generalized Langevin equations (GLE),
\begin{equation}\label{eq: gle}
 \left\{
  \begin{aligned}
  	\dot{x}=& v, \\
	\dot{v} = &f(x) - \gamma v - \int_0^t \theta(t-\tau) v(\tau) d\tau + \xi(t).
   \end{aligned}\right.
\end{equation}
The mass has been set to unity. The random noise is assumed to be a mean-zero Gaussian process, with time correlation given by,
\begin{equation}\label{eq: fdt0}
  \Big\langle \xi(t) \xi(t')^T \Big\rangle = 2k_BT\gamma \delta(t-t') + k_B T \theta(t-t').
\end{equation}
This is a general form of the fluctuation-dissipation theorem (FDT) \cite{Kubo66}. Equation of this form has been derived in our previous work \cite{MaLiLiu16}, where the full model is the Langevin dynamics.

We first assume that $\theta$ and $\gamma$ are symmetric, matrix-valued functions. 
Our goal is simply solving the second equation, and then making a direct substitution into the first equation to eliminate $v(t)$ entirely.
For this purpose, we define a matrix-valued function $\chi,$ which satisfies the following equation,
\begin{equation}\label{eq: def-chi}
 \dot{\chi}= - \gamma \chi - \int_0^t \theta(t-\tau) \chi(\tau) d\tau, \quad \chi(0)=I.
\end{equation}

Through Laplace transform, one can show that the function $\chi$ is symmetric. Now $\chi$ can be viewed as the fundamental function associated with the GLE model \eqref{eq: gle}. In particular, we have,
\begin{equation}
  v(t)= \chi(t) v(0) + \int_0^t \chi(t-\tau) \xi(\tau) d\tau + \int_0^t \chi(t-\tau) f\big(x(\tau)\big) d\tau.
\end{equation}

Let us first look at the first two terms, denoted by $w(t)$
\begin{equation}\label{eq: fdt1}
 w(t)= \chi(t) v(0) + \int_0^t \chi(t-\tau) \xi(\tau) d\tau.
\end{equation}
We assume that $v(0)$ is sampled from its equilibrium. Namely it is Gaussian with mean zero and covariance $k_B T I.$ Then the following result can be established using the FDT \eqref{eq: fdt0}: $w(t)$ is a stationary Gaussian process with mean zero and time correlation given by,
\begin{equation}
\Big\langle w(t) w(t')^T \Big\rangle= k_B T \chi(t-t'),\; \forall t\ge t'. \label{wcor}
\end{equation}
The proof is postponed to the appendix due to some non-trivial calculations.

\medskip

With this result, we arrive at a Brownian dynamics {\it with memory},
\begin{equation}\label{eq: bd-mem}
  \dot{x}= \int_0^t \chi(t-\tau) f\big(x(\tau)\big) d\tau + w(t).
\end{equation}
Together with \eqref{eq: fdt1}, this forms a closed set of equations, since a stationary Gaussian process is uniquely determined by its time correlation \cite{Doob44}. 

To arrive at a conventional BD model \eqref{eq: bd}, we may approximate the kernel function by,
\begin{equation}\label{eq: chi-inf}
 \chi(t) \approx 2\chi_\infty \delta(t),\quad  \chi_\infty = \int_0^{+\infty} \chi(t) dt,
\end{equation}
and we obtain,
\begin{equation}\label{eq: bd1}
 \dot{x}(t) = \chi_\infty f(x(t)) + w(t).
\end{equation}

In order to preserve the FDT \eqref{wcor}, the noise now has to be approximated by a white noise, with covariance given by,
\begin{equation}
 \Big\langle w(t) w(t')^T \Big\rangle= 2k_B T \chi_\infty \delta(t-t').
\end{equation}
This is exactly the standard Brownian dynamics model (\ref{eq: bd}). The coefficients are determined from the formula \eqref{eq: chi-inf}, which is clearly of linear-response type. {\blue Namely, we have,
\[\chi_\infty =\frac{1}{2k_B T } \int_0^{+\infty} \Big\langle w(t) w(0)^T \Big\rangle dt.\]}
 In the next section, we will examine the validity of this approximation.

\section{A Case Study}
To test the approximation \eqref{eq: bd1} presented in the previous section, we consider the GLE model derived by Adelman and Doll \cite{AdDo74,AdDo76}, in which a gas molecule near a solid surface was considered. We write the position of the gas molecule as $y(t)$ and the position of the solid atom at the surface as $x(t).$ Starting from a one-dimensional chain of solid atoms, Adelman and Doll \cite{AdDo74} used Laplace transform, together with a spatial reduction procedure, and derived a GLE for $x(t)$ and $y(t)$ with the remaining solid atoms eliminated. The GLE can be written as,
\begin{equation}\label{eq: Doll}
\left\{
 \begin{aligned}
   M \ddot{y}=& - \phi'(x-y),\\
   \ddot{x} = & \phi'(x-y) -\gamma \dot{x} - \int_0^t \theta(t-\tau) \dot{x}(\tau) d\tau + R(t).
 \end{aligned}\right.
\end{equation}
Here $M$ is the mass of the gas molecule, and as in the original work \cite{AdDo74}, we set the mass of the solid atom to unity. We write the memory term as a friction term with $\dot{x}$ involved. This way, the FDT is directly expressed as,
\begin{equation}
 \big\langle R(t) R({\blue t'})\big\rangle=k_BT \theta(t-t')+ 2k_B T \gamma \delta(t-t').
\end{equation}
We also add the damping term with coefficient $\gamma$ to offer a slightly more general setup.

When the one-dimensional chain of solid atoms is infinite, an explicit formula for the kernel function can be found via Laplace transform, and it is given by,
\begin{equation}\label{eq: theta}
 \theta(t) = \frac{\omega_0}{ t} J_1(2 \omega_0 t),
\end{equation}
where $J_1$ is the Bessel function of the first kind, and $\omega_0=\sqrt{K}$ with $K$ being the spring constant of the atom chain. This model has been considered by \cite{LiE07,LiE06,liu2006nano} to study boundary conditions for molecular dynamics models. Although this explicit formula might be already known to others,  we have not found earlier references to it. This kernel function has very slow decay rate, which makes it a highly non-trivial test problem. In light of the natural time scaling in \eqref{eq: theta}, we will consider $\omega_0$ as the natural frequency associated with the kernel function, and it represents the corresponding time scale. To make this model more explicit, we choose the Morse potential for the interaction of the gas/solid molecules at the interface,
\begin{equation}
 \varphi(u)= (e^{-au}-1)^2. 
\end{equation}
$a$ is chosen as 1, so the vibration frequency associated with the force $f(x)$ is $\omega=\sqrt{2}.$  To further simplify the model, we fix $y=0,$ so that only the second equation in (\ref{eq: Doll}) needs to be solved.

The Laplace transform of the kernel function is given by,
\begin{equation}\label{eq: theta-lap}
 \Theta(s)=  \frac{\sqrt{s^2+4\omega^2_0}}2 -\frac{s}2.
\end{equation}
This will be needed in the approximation schemes presented here and the ones in the next section. In particular, we obtain,
\begin{equation}\label{eq: chi-inf-1d}
 \chi_{\infty} = \frac{1}{\gamma +\omega_0}.
\end{equation}

\begin{figure}[htbp]
\begin{center}
\includegraphics[scale=0.5]{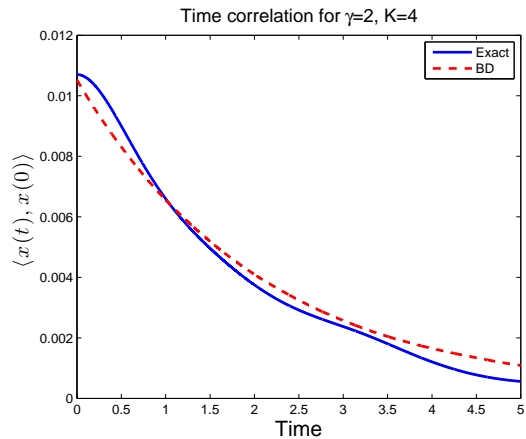}
\caption{Comparison of the time correlation functions for the full model and the BD model \eqref{eq: bd1}. $\gamma=2, K=4$.}
\label{fig: BD-0}
\end{center}
\end{figure}

To examine the accuracy, we first run the full dynamics with 8192 solid atoms and generate the trajectory of the solid atom at the gas/solid interface. We use the stochastic velocity-Verlet method \cite{van1982algorithms} to simulate the Langevin dynamics model.  The time correlation function $\bl x(t), x(0) \br$ as a dynamic property is regarded  as the exact result.  As comparison, we solve the BD model \eqref{eq: bd} with the effective damping coefficients given by \eqref{eq: chi-inf-1d} using the Euler-Maruyama method. We first choose $\gamma=2,$ and $K=4$, which implies $\omega_0=2$. Figure \ref{fig: BD-0} displayed the two time correlation functions. One can see that the results agree quite well, indicating that in this regime, the approximation by the BD model \eqref{eq: bd} is good. 

\begin{figure}[htbp]
\begin{center}
\includegraphics[scale=0.5]{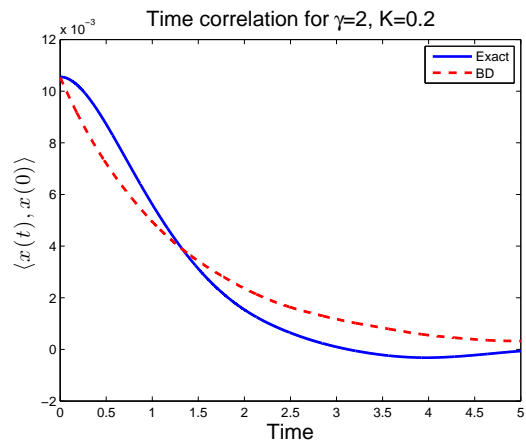}
\caption{Comparison of the time correlation functions for the full model and the BD model \eqref{eq: bd1}. $\gamma=2, K=0.2$.}
\label{fig: BD-4}
\end{center}
\end{figure}

In the next test, we set $K=0.2$, which gives a much lower frequency. The results are shown in Figure \ref{fig: BD-4}, from which we observe discrepancies between the results. In this case, the approximation leads to much bigger error. 

\begin{figure}[htb]
\begin{center}
\includegraphics[scale=0.4]{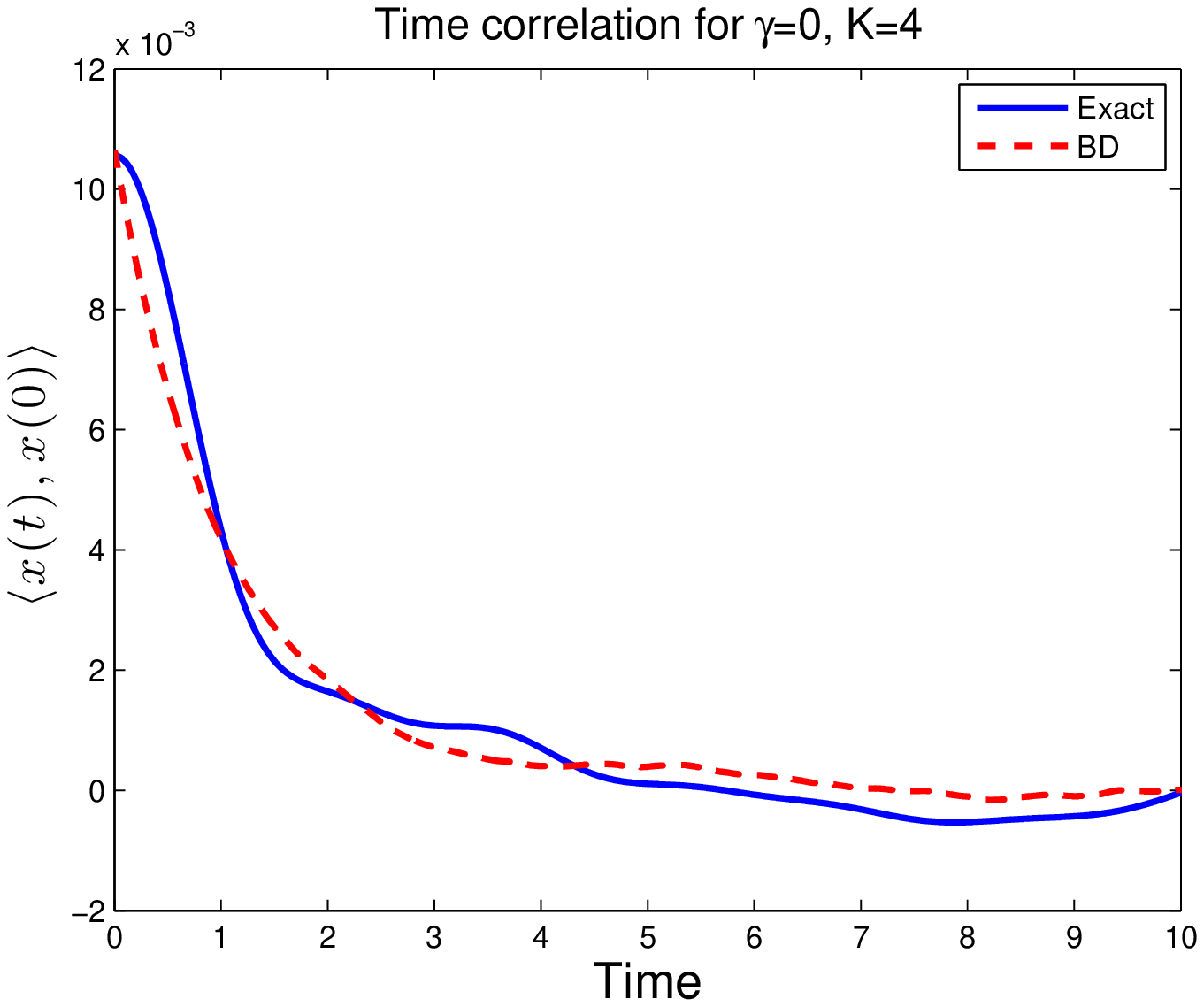}
\includegraphics[scale=0.4]{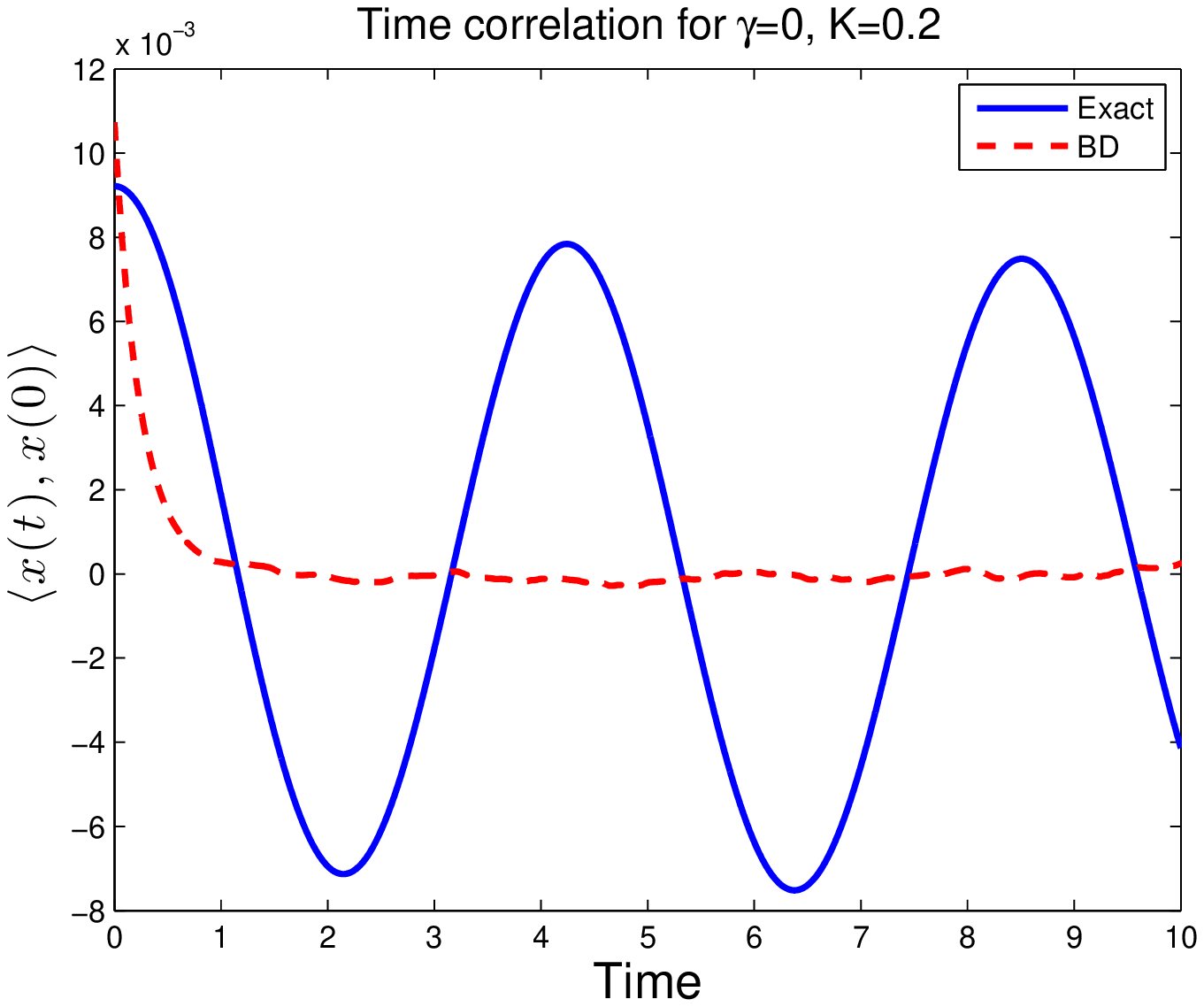}
\caption{Comparison of the time correlation functions for $\gamma=0,$  $K=4$ (left) and $K=0.2$ (right).}
\label{fig: gam0}
\end{center}
\end{figure}

Next, we repeat these experiments with $\gamma=0.$ In this case, the damping mechanism, e.g., $\chi_\infty$, can be directly attributed to the  friction kernel $\theta.$ When $K=4,$ it turns out that BD is still a reasonable approximation. However, when $K=0.2,$ the time correlation exhibits pronounced oscillations, which is not captured by the BD model, as shown in Figure \ref{fig: gam0}.

\section{Embedded Brownian Dynamics}

It has been observed from the previous section that in the regime where $\omega_0 \ll \omega$, the Brownian dynamics might not be a good approximation
to the original dynamics. In this case, we propose to embed the reduced model \eqref{eq: bd-mem} in an extended system, where the memory is removed and the general Gaussian noise is generated from the additional equations. 

More specifically, the approximation is done in terms of  the Laplace transform of the kernel function, given by,
\begin{equation}
 X(s)= \Big[sI + \gamma + \Theta(s)\Big]^{-1}.
\end{equation}
We seek a rational approximation, in the form of $R(s)= Q(s)^{-1}P(s),$ where $P$ and $Q$ are matrix polynomials, given by,
\begin{equation}
\begin{aligned}
P(s)=& A_0 s^{n-1} + A_1 s^{n-1} + \cdots + A_n,\\
Q(s)=& s^n I - B_0 s^{n-1} - B_1 s^{n-1} + \cdots + B_n. 
\end{aligned}
\end{equation}
Here $n$ indicates the level of approximation.

For the zeroth level approximation, we have a constant function and we choose to match the limit as $s\to 0.$ We obtain,
\begin{equation}
  \lim_{s\to 0} X= \Big[\gamma + M_\infty\Big]^{-1},
\end{equation}
with $M_{\infty}=\lim_{s\to 0} \Theta$. This limit is defined as $\chi_\infty$. Namely, we choose $R(s)=\chi_\infty.$ With this constant approximation in the Laplace domain, we have a corresponding approximation in the time domain,
which gives exactly the BD model \eqref{eq: bd1}.

We now turn to higher-order rational approximations.
For $n=1$, we write $R_{1,1}(\lambda)=[1-\lambda B]^{-1} A\lambda$, where $\lambda=1/s.$
Similarly we rewrite $X$ as a function of $\lambda$,
\begin{equation}
 X(\lambda)= \Big[I + \lambda \gamma + \lambda \Theta(\lambda)\Big]^{-1} \lambda,
\end{equation}
and approximate this Laplace transform by $R_{1,1}$.

From the definition of the Laplace transform, it is clear that $X, R \to 0$, as $\lambda \to 0.$  To determine the matrices $A$ and $B$, we impose the following matching conditions:
\begin{align}\label{eq: matching-cond}
\begin{split}
&\lim_{\lambda \to +\infty} X(\lambda)=\lim_{\lambda \to +\infty} R_{1,1}(\lambda),\\
&\lim_{\lambda \to 0_+} \frac{d}{d\lambda}X(\lambda)=\lim_{\lambda \to 0_+} \frac{d}{d\lambda}R_{1,1}(\lambda).
\end{split}
\end{align}

Direct calculations yield, 
\begin{equation}\label{eq: n=1}
A=I, \quad B=-(\gamma+M_{\infty}).
\end{equation}
This is equivalent to approximating $\chi(t)$ as matrix-exponential function,
\begin{align}\label{eq: chi1}
\chi(t)\approx \chi_I(t) \overset{\text{def}}{=} e^{Bt}.
\end{align}

The rational approximation in the Laplace domain makes it possible to eliminate the memory by introducing auxiliary variables. We set $z=\int_0^t \chi_I(t-\tau)f(x(\tau))d\tau$ in (\ref{eq: bd}) as the memory function, then we  write an ODE for $z$:
\begin{align}\label{eq: z0}
&\dot{z}=B z + f, \quad z(0)=0.
\end{align}

We now show that by adding a white noise in this equation, i.e., we have,
\begin{equation}\label{eq: 1st-order}
\left\{
\begin{aligned}
&\dot{x}=z,\\
&\dot{z}=Bz+f(x)+\xi(t),\\
\end{aligned}\right.
\end{equation}
the Gaussian noise $w(t)$ in the exact reduced model \eqref{eq: bd-mem} can be approximated so that the correlation function is {\it exactly} 
consistent with the  approximation by $\chi_I(t).$

To prove this statement, let $\Sigma$ be the covariance the noise, i.e., 
$$
\langle \xi(t)\xi(t')^T\rangle = \Sigma\delta(t-t').
$$
Then by solving $z$ analytically, we have the expression:
\begin{align*}
z(t) = e^{Bt} z(0) +\int_0^t e^{B(t-\tau)} \xi(\tau) d\tau+\int_0^t e^{B(t-\tau)} f(x(\tau)) d\tau.
\end{align*}
Notice the third term of $z(t)$ is exactly the approximation of the memory term with the memory kernel given by $\chi_I(t)$. 
Let the covariance of $z(0)$ be $Q=k_B T I$. Then the first two terms form a stationary Gaussian process, provided that the covariance $\Sigma$ satisfies the Lyapunov equation \cite{risken1984fokker}:
\begin{align}
- \Sigma=BQ+QB^T.
\end{align}
Furthermore, the time correlation of the stationary noise is given by $e^{B(t-t')}Q$, which is exactly $\chi_I(t)$, therefore  the FDT as in 
\eqref{wcor} is satisfied.

\bigskip

It is clear that this procedure can be extended to higher order. For example, to construct the approximation for $n=2$, we set,
\begin{equation}\label{eq: R22}
R_{2,2}=[I-\lambda B_0-\lambda^2 B_1]^{-1}[\lambda A_0 + \lambda^2A_1].
\end{equation}
To determine the coefficient matrices, we impose the following four conditions,
\begin{equation}\label{eq: conds-2}
\begin{aligned}
&\lim_{\lambda \to +\infty} X(\lambda)=\lim_{\lambda \to +\infty} R_{1,1}(\lambda),\\
&\lim_{\lambda \to 0_+} \frac{d^\ell}{d\lambda^\ell}X(\lambda)=\lim_{\lambda \to 0_+} \frac{d^\ell}{d\lambda^\ell}R_{1,1}(\lambda), \quad \ell=1,2,3.
\end{aligned}
\end{equation}
The explicit form of the solutions for the coefficients $A_0, A_1, B_0,$ and $B_1$ can be found in the appendix.

We let the corresponding approximate kernel function be $\chi_{\textrm{II}}(t),$ which is the inverse Laplace transform of $R_{2,2}.$ Similar to the previous approximation, the memory term can be embedded into an extended system with auxiliary variables $z$ and $z_1,$ together with the white noises $\xi_1$ and $\xi.$ {\blue The extended system of equations read,}
\begin{equation}\left\{
\begin{aligned}\label{eq: n=2}
&\dot{x}=z,\\
&\dot{z}_1=B_1 z +A_1 f + \xi_1,\\
&\dot{z}=z_1 + B_0 z + A_0 f + \xi.
\end{aligned}\right.
\end{equation}

Through Laplace transform and direct substitutions, one can easily check that the memory function $\chi$ is approximated by $\chi_{\textrm{II}}(t)$, with the inverse Laplace transform given by \eqref{eq: R22}. The random noise $w(t)$ is also approximated by a Gaussian noise that satisfies the FDT \eqref{wcor} with proper choices of the initial conditions for $z(t)$ and $z_1(t)$ and the covariance of the noises $\xi_1$ and $\xi$. {\blue The details can be found in Appendix C.}

\begin{figure}[htbp]
\begin{center}
\includegraphics[scale=0.4]{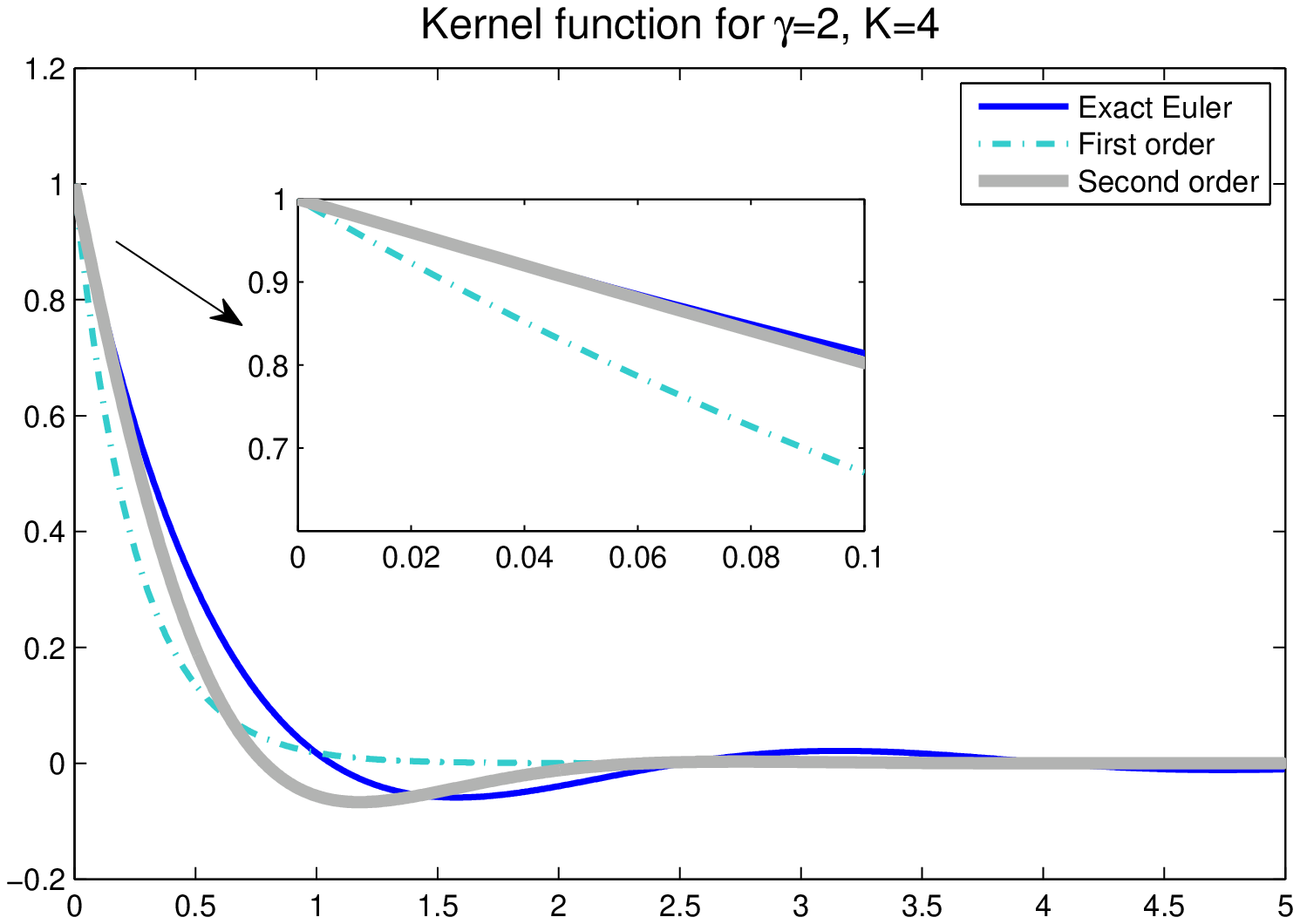}
\includegraphics[scale=0.4]{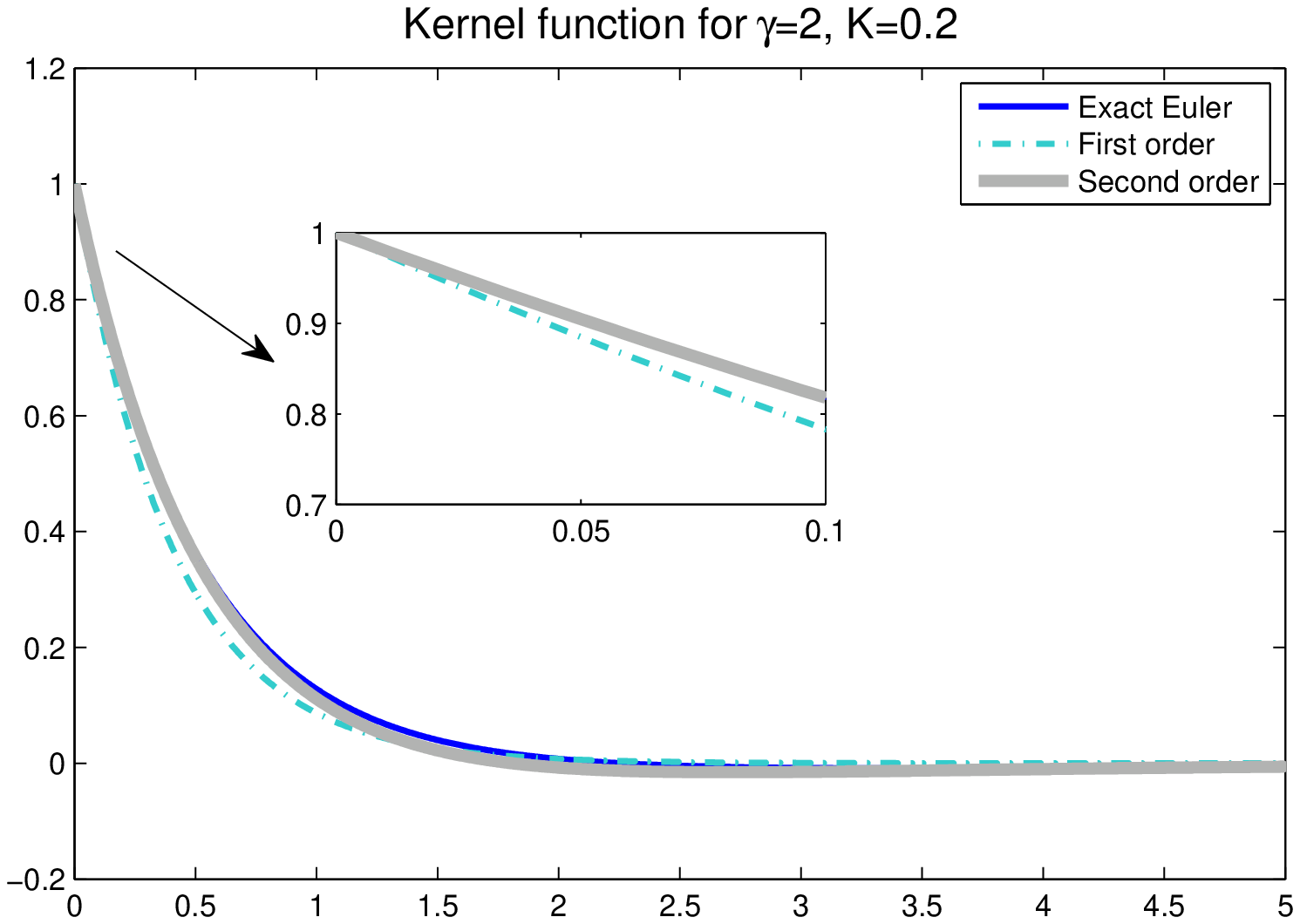}
\caption{Comparison of the kernel functions for the exact result,  the first order rational approximation  \eqref{eq: chi1}, and the second order approximation \eqref{eq: chi2}. $\gamma=2$, $K=4$ (left), $K=0.2$ (right).}
\label{fig: ker}
\end{center}

\end{figure}

Now we turn to the test problem and compare the results from these models. In Figure \ref{fig: ker}, we present a comparison between different orders of approximations to the exact kernel function for $\gamma=2$. The left panel is for $K=4$ while the right panel is for $K=0.2$. We did not plot the zeroth order approximation {\blue (the BD model \eqref{eq: bd1}) } since it is a delta function. We want to comment on how we produced the exact result. While the explicit expression of $\theta(t)$ is stated in (\ref{eq: theta}), the exact form of $\chi(t)$ is however nontrivial, which involves integral of highly oscillatory Bessel functions. We therefore adopt a numerical Euler method to invert Laplace transform \cite{abate2006unified} and name it as Exact Euler. {\blue This exact result is tested against high-order numerical quadrature to ensure its accuracy.}

\begin{figure}[htbp]
\begin{center}
\includegraphics[scale=0.45]{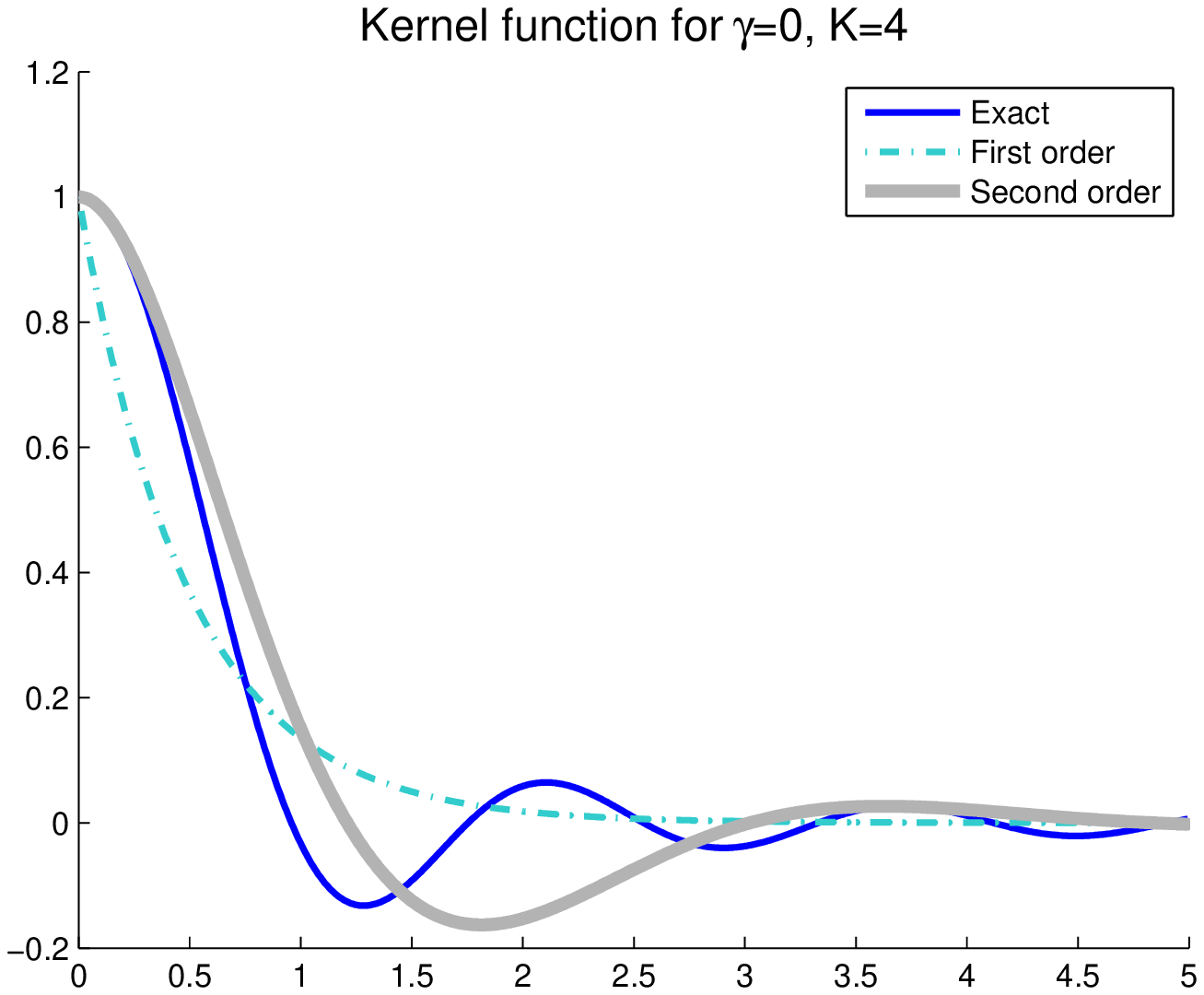}
\includegraphics[scale=0.45]{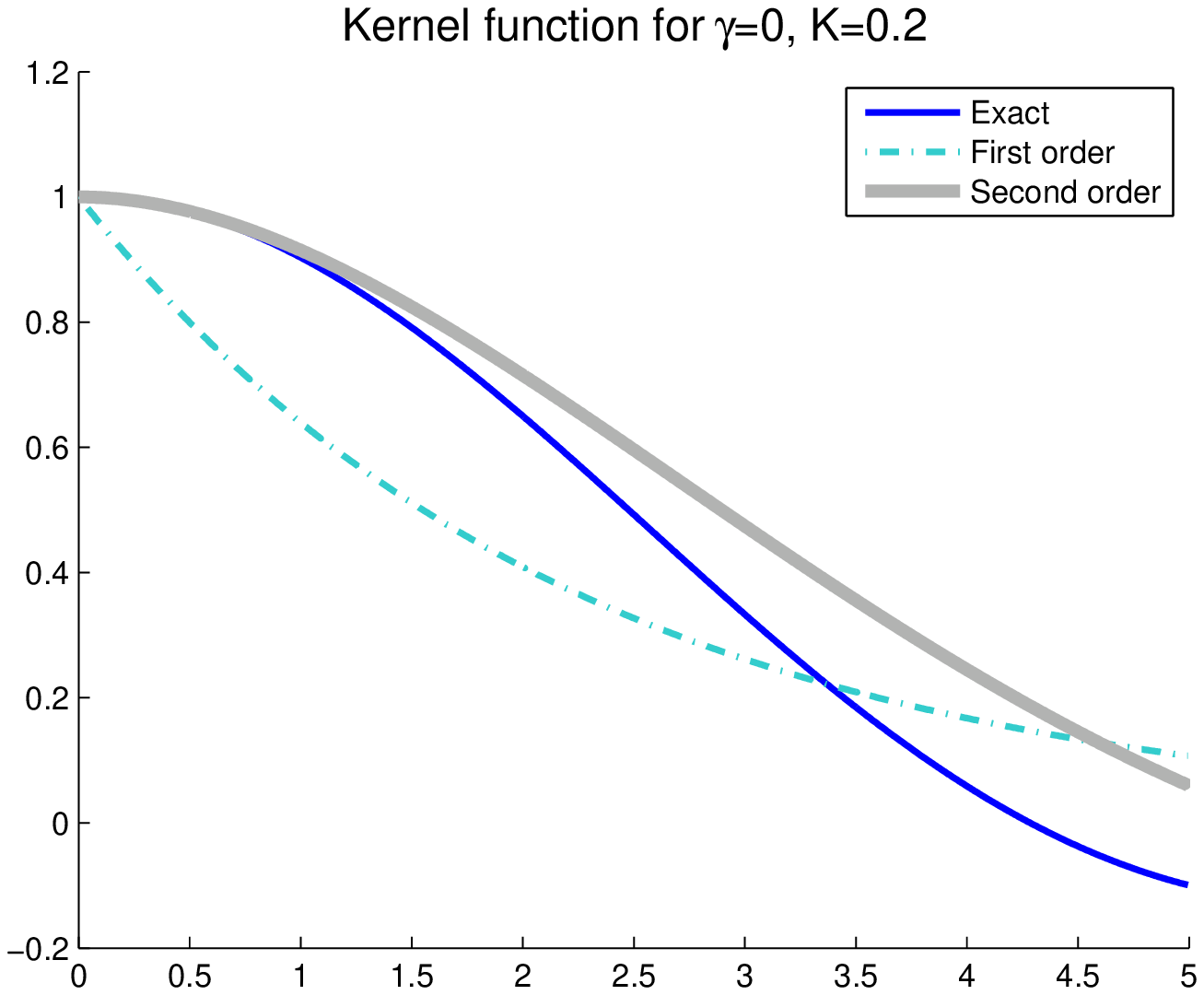}
\caption{Comparison of the kernel functions.  $\gamma=0$, $K=4$ (left), $K=0.2$ (right).}
\label{fig: kergam0}
\end{center}
\end{figure}

One can observe from the graph in Figure \ref{fig: ker} that for smaller frequency in the memory kernel $\theta(t)$, the kernel function $\chi(t)$ is more flat, {\blue as shown in the graph on the right}, and we have improved accuracy from the rational approximations $R_{1,1}$.  At the same time, it is obvious that with higher order of approximation $R_{2,2}$, we achieve better results. Next we repeat test for $\gamma=0$. In this case, $\chi$ can be explicitly expressed as
$$
\chi(t)=\frac{ J_1(2\omega_0 t)}{t\omega_0}.
$$ 
Results are shown in Figure \ref{fig: kergam0},  and {\blue it is clear that for this case, much larger error is introduced in the first order  approximation. Though  the second order approximation offers better agreement, there is still obvious discrepancy. Therefore we carried out a third order approximation in a similar algorithm, and clearly observe the tendency of convergence to the exact kernel function. Similar to the situation for $\gamma=4$, for smaller frequency in $\theta(t)$, the kernel function is more flat.}


\begin{figure}[htbp]
\begin{center}
\includegraphics[scale=0.5]{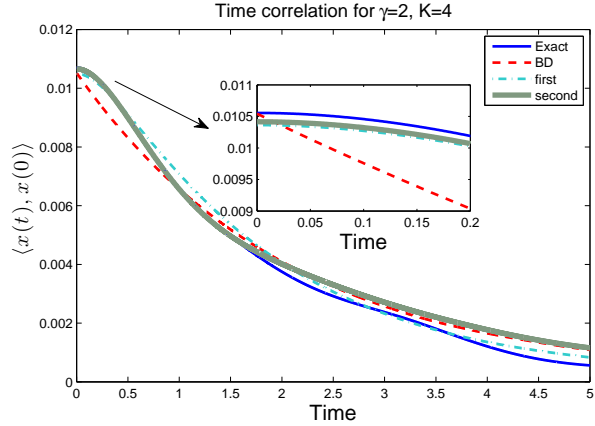}
\caption{Comparison of the time correlation functions for the full model, the BD model \eqref{eq: bd}, the first order rational approximation  \eqref{eq: 1st-order}, and the second order approximation \eqref{eq: n=2}. $\gamma=2, K=4$.}
\label{fig: all0}
\end{center}
\end{figure}

Next, we compare the time correlation of all the cases mentioned. For $\gamma=2$, $K=4,$ the results are collected in Figure \ref{fig: all0}. In this case, all the approximations, including the BD model \eqref{eq: bd1}, the first order rational approximation  \eqref{eq: 1st-order}, and the second order approximation \eqref{eq: n=2}, yield reasonable results. The second order method offers the best approximation over short time scale, {\blue this is consistent with increased accuracy near $t=0$ in the approximations}.   

\begin{figure}[htbp]
\begin{center}
\includegraphics[scale=0.5]{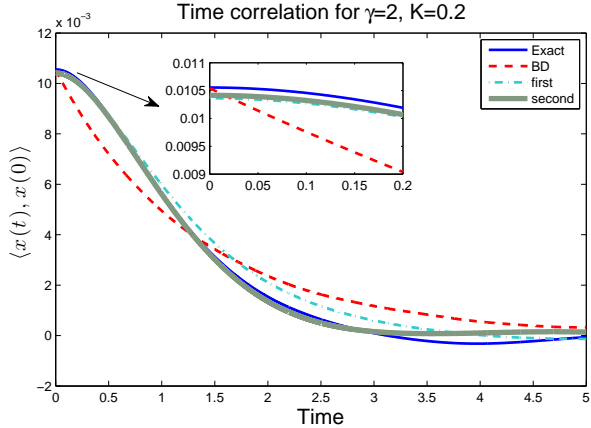}
\caption{Comparison of the time correlation functions. $\gamma=2, K=0.2$.}
\label{fig: all4}
\end{center}
\end{figure}

Meanwhile, for a lower frequency from the kernel, where $K=0.2,$ the results are shown in Figure \ref{fig: all4}. We observe that as the order of approximation increases, the accuracy is significantly improved. The insets in both Figure \ref{fig: all0} and \ref{fig: all4} showed that, for higher order of approximations, derivatives near $t=0$ are better predicted.  

\begin{figure}[htbp]
\begin{center}
\includegraphics[scale=0.45]{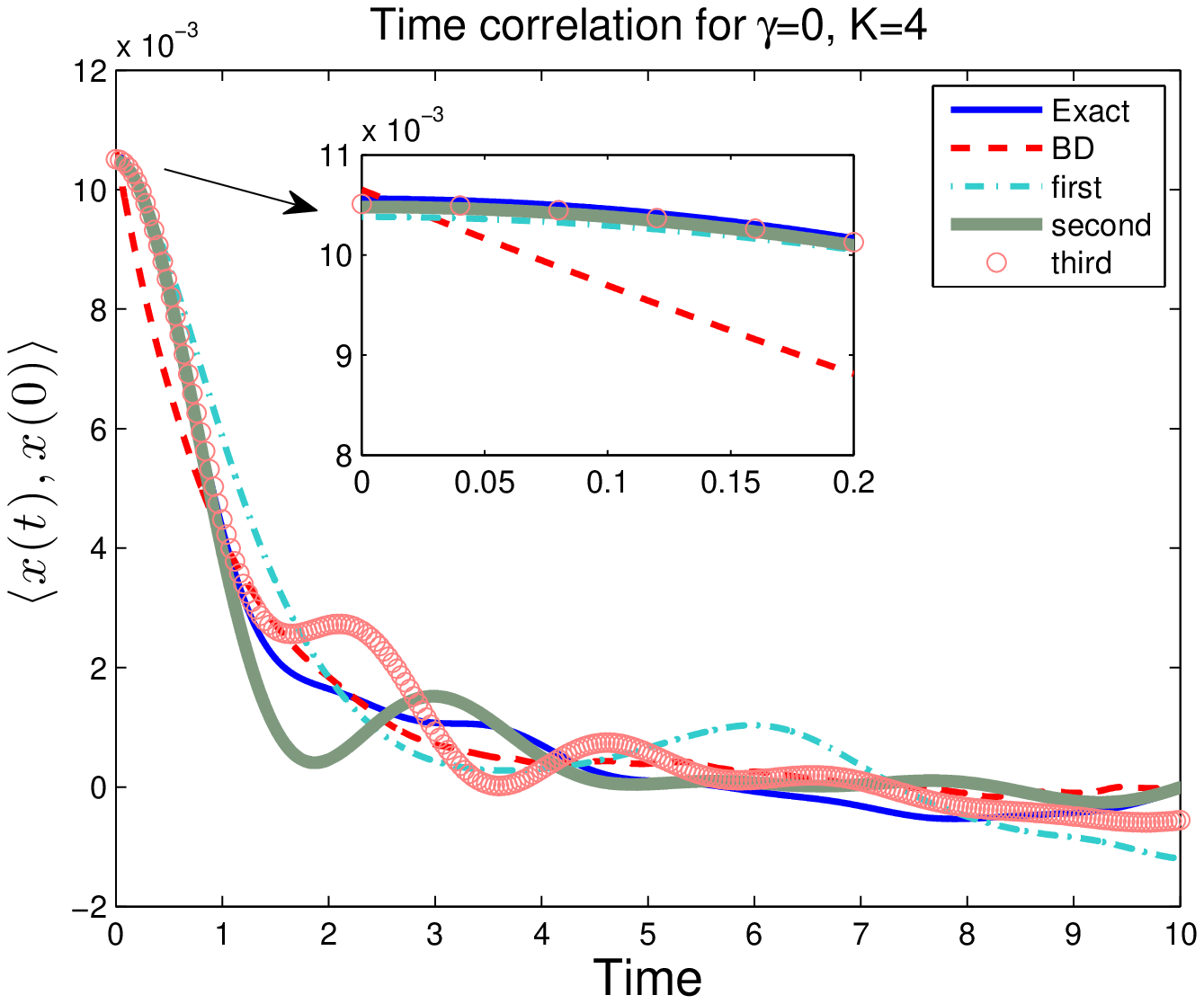}
\includegraphics[scale=0.45]{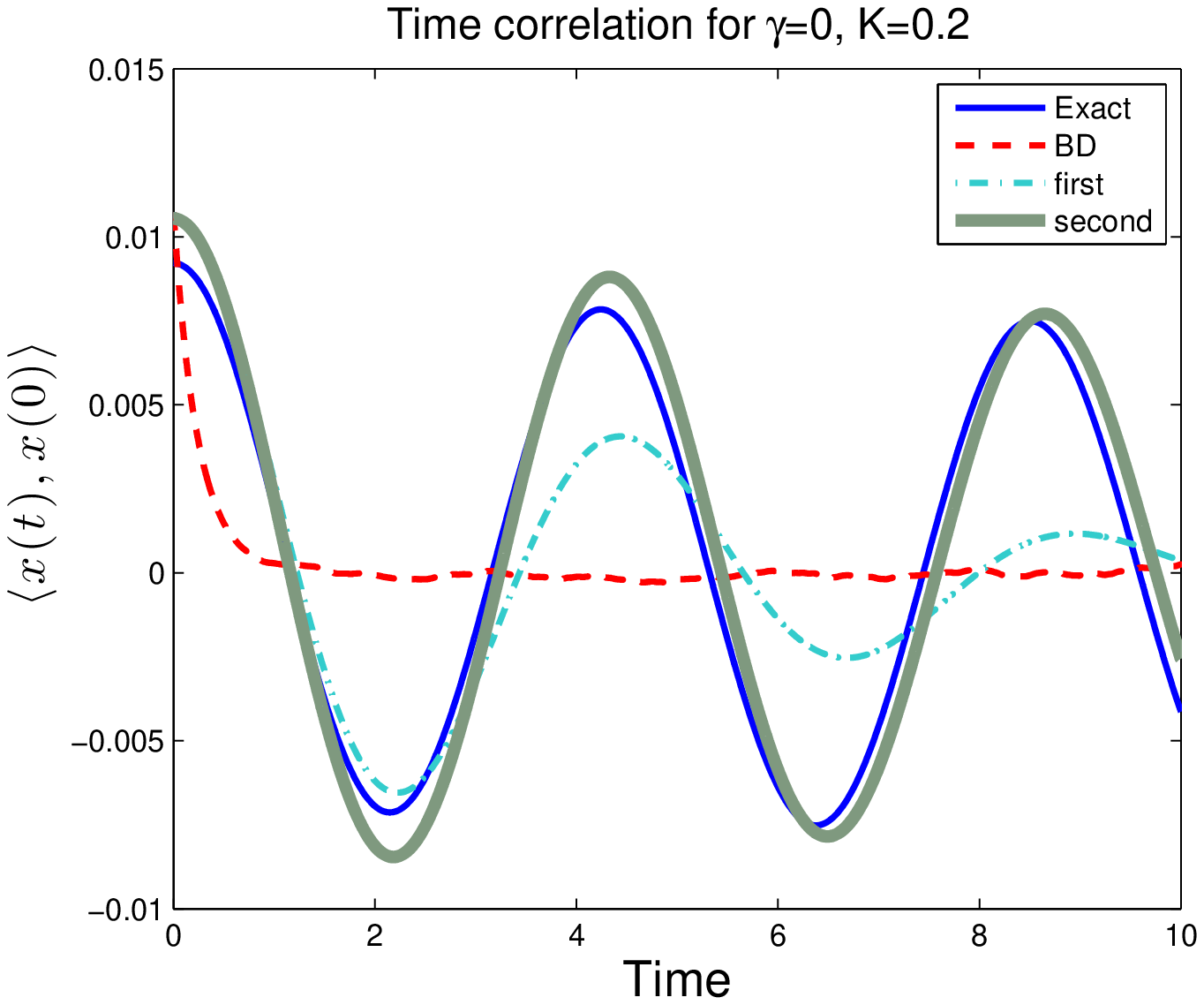}
\caption{Comparison of the time correlation functions among the approximate models for $\gamma=0,$  $K=4$ (left) and $K=0.2$ (right).}
\label{fig: gam0_all}
\end{center}
\end{figure}

We  also  repeated these experiments with $\gamma=0.$  When $K=4,$ all the approximations give reasonable predictions. {\blue One can clearly see the improvement near zero as approximation order gets higher, due to more conditions are matched at zero.} However, when $K=0.2,$ the large oscillations in the time correlation are not captured by the BD model as shown in Figure \ref{fig: gam0_all}. The first order method does predict oscillations, but with the wrong magnitude. The second order method provides the better approximation, {\blue and third order is definitely the best among all}. 

\section{Summary and Discussions}
In this paper we considered the reduction of the generalized Langevin equations to a Brownian dynamics with memory (nonlocal) and correlated noise, then to the standard Brownian dynamics as well as higher order approximations.  
Due to the considerable advantage of the Brownian dynamics model over Langevin or generalized Langevin dynamics in allowing much larger time steps and longer time scale, it is important to understand the accuracy of such approximations. Our findings gave insight into when to use the simplest BD model and how to reduce the error by using high order approximations.  

In our numerical experiments, we chose the memory kernel of  the GLE model derived by Adelman and Doll \cite{AdDo74}. The corresponding memory kernel exhibits very slow decay, which makes it a good candidate to understand the modeling error. Our interest was in finding the appropriate regime where the Brownian dynamics is a good approximation for the GLE model. In particular, numerous experiments with  various   frequencies ($\omega_0$) of the memory kernel were conducted, while fixing the vibration frequency $\omega$ associated with the potential of mean forces. Our observations have been that the Brownian dynamics is a reasonably good approximation in  cases when $\omega_0$ is relatively large. 

Meanwhile, the Brownian dynamics model will generate large error in the case of low frequency in the memory kernel. We derived high order approximations, by embedding the nonlocal Brownian dynamics model into an extended, {\it local}  stochastic dynamics. We have shown that by correctly choosing the covariance for the additive noises and the auxiliary variables, the second fluctuation-dissipation theorem are exactly satisfied. Our numerical results showed improved accuracy as the order of the approximation is increased.


At the same time, the current approach for constructing high order approximations can be extended/improved. For example, in constructing our models, we only used interpolations of the Laplace transform at $\lambda=0$ and $\lambda=+\infty$. It is possible to introduce more interpolation points to better resolve the intermediate rescales.
 
\section{Acknowledgment} 
This research was supported by NSF under grant DMS-1412005, DMS-1216938 and DMS-1619661. 

\section{Appendix}
\subsection{Stationary Gaussian process $w(t)$}

In this section, we will prove the statement made in section \ref{eq: sec-deri}. Namely, we will  show that $w(t)$ is indeed a stationary Gaussian process satisfying the stated property.  More specifically,  we want to show  that 
$\Big\langle w(t) w(t')^T\Big \rangle $ is a function of $(t-t')$, which  is equivalent to showing 
$$
\frac{d}{ds} \Big\langle w(t+s) w(t'+s)^T\Big \rangle = 0,
$$
for any $s.$ Consequently,  by taking $s=-t'$,
$$
\Big\langle w(t) w(t')^T\Big \rangle = \Big\langle w(t-t') w(0)^T\Big \rangle,
$$
which shows the stationarity of the random process.

Recall that 
\begin{equation*}
 w(t)= \chi(t) v(0) + \int_0^t \chi(t-\tau) \xi(\tau) d\tau,
\end{equation*}
and 
\begin{equation*}
  \Big\langle \xi(t) \xi(t')^T \Big\rangle = 2k_BT\gamma \delta(t-t') + k_B T \theta(t-t'),
\end{equation*}
together with the property that $\chi(t)$ is symmetric. For $t\geq t'$, the correlation can be written as 
\begin{align*}
&\Big \langle w(t+s) w(t'+s)^T \Big \rangle \\
=& \chi(t+s) \Big\langle v(0) v(0)^T\Big \rangle \chi(t'+s) + 
\int_0^{t+s} \int_0^{t'+s} \chi(t+s-\tau) \Big \langle \xi(\tau) \xi(\tau')^T\Big\rangle \chi(t'+s-\tau') d\tau' d\tau\\
=&k_BT\chi(t+s)\chi(t'+s)\\
&+ 2k_BT\int_0^{t'+s} \chi(t+s-\tau') \Gamma \chi(t'+s-\tau') d\tau'\\
 &+ 
k_BT\int_0^{t+s}\int_0^{t'+s} \chi(t+s-\tau)\theta(\tau-\tau') \chi(t'+s-\tau') d\tau' d\tau.
\end{align*}
We denote these terms  as term I, \textrm{II}, \textrm{III} respectively.

Taking derivatives with respect to $s$ to each term, we find:
For Term I, we have:
\begin{align*}
&\frac{d}{ds}\chi(t+s)\chi(t'+s) = \dot{\chi}(t+s) \chi(t'+s) + \chi(t+s)\dot{\chi}(t'+s)\\
=& \left(- \chi(t+s) \Gamma - \int_0^{t+s} \chi(\tau) \theta(t+s-\tau) d\tau \right) \chi(t'+s)
+ \chi(t+s) \left(-\Gamma \chi(t'+s) - \int_0^{t'+s} \theta(t'+s-\tau)\chi(\tau) d\tau\right)\\
=&-2\chi(t+s)\Gamma\chi(t'+s) - \int_0^{t+s} \chi(\tau) \theta(t+s-\tau)\chi(t'+s) d\tau
- \int_0^{t'+s} \chi(t+s) \theta(t'+s-\tau) \chi(\tau) d\tau,
\end{align*}

for Term \textrm{II}, we have:
\begin{align*}
&2\frac{d}{ds} \int_0^{t'+s} \chi(t+s-\tau') \Gamma \chi(t'+s-\tau') d\tau' = 2\chi(t+s-(t'+s))\Gamma\chi(t'+s-(t'+s))\\
&+ 2\int_0^{t'+s} \dot{\chi}(t+s-\tau')\Gamma \chi(t'+s-\tau') d\tau+
 2\int_0^{t'+s} {\chi}(t+s-\tau')\Gamma \dot{\chi}(t'+s-\tau') d\tau'\\
 =& 2\chi(t-t')\Gamma\chi(0) -2 \int_0^{t'+s} \frac{d}{d\tau'} \left(\chi(t+s-\tau')\Gamma \chi(t'+s-\tau')\right) d\tau'\\
 =& 2\chi(t-t')\Gamma\chi(0) - 2 \chi(t+s-\tau')\Gamma \chi(t'+s-\tau')\Big|_{\tau'=0}^{t'+s}\\
 =& 2\chi(t+s)\Gamma\chi(t'+s),
\end{align*}

and for Term \textrm{III}:
\begin{align*}
&\frac{d}{ds}\int_0^{t+s}\int_0^{t'+s} \chi(t+s-\tau)\theta(\tau-\tau') \chi(t'+s-\tau') d\tau' d\tau\\
\xlongequal{\text{$u=t'+s-\tau'$}}&
\frac{d}{ds}\int_0^{t+s}\int_0^{t'+s} \chi(t+s-\tau)\theta(\tau-(t'+s-u)) \chi(u) d u d\tau\\
\xlongequal{\text{$v=t+s-\tau$}}&
\frac{d}{ds}\int_0^{t+s}\int_0^{t'+s} \chi(v)\theta(t-t'-v+u) \chi(u) d u dv\\
=&\int_0^{t'+s}\chi(t+s)\theta(t-t'-(t+s)+u) \chi(u) d u + \int_0^{t+s} \chi(v)\theta(t-t'-v+(t'+s)) \chi(t'+s) dv\\
=&\int_0^{t'+s}\chi(t+s)\theta(t'+s-\tau) \chi(\tau) d \tau+ \int_0^{t+s} \chi(\tau)\theta(t+s-\tau)) \chi(t'+s) d\tau.
\end{align*}
Notice that the last equation used the property that $\theta$ is an even function. It is clear that the summation of all the three terms gives zero, which shows the stationarity. 

Now we have for $t\geq t'$,
\begin{align*}
&\Big\langle w(t) w(t')^T\Big \rangle = \Big\langle w(t-t') w(0)^T\Big \rangle\\
=& \Big\langle\left (\chi(t-t')v(0)+\int_0^{t-t'} \chi(t-t'-\tau)\xi(\tau) d\tau \right)v(0)^T \Big \rangle = k_BT \chi(t-t').
\end{align*}
This is precisely the stated result. \qed

\subsection{The derivation of the second order system \eqref{eq: n=2}}
We approximate the Laplace transform of $\chi$ in the following form:
\begin{align*}
(s^2-sB_0-B_1) X= sA_0 + A_1,
\end{align*}
this is equivalent as
\begin{align*}
&s X_1 - B_1 X = A_1,\\
& X_1 = sX-B_0X - A_0.
\end{align*}
Taking the inverse of Laplace transform, we arrive at the differential equations:
\begin{align*}
&\dot{\chi_1} = B_1 \chi, \quad \chi_1(0) = A_1,\\
&\dot{\chi} = B_0 \chi + \chi_1, \quad \chi(0)=A_0.
\end{align*}
{\blue In this calculation, we have used the matching conditions at $\lambda=0_+$ \eqref{eq: conds-2}. Therefore, incorporating the values of $X(\lambda)$ at $\lambda=0_+$ seems to be important to convert the model to the time domain.}

The solution $\chi$ will be denoted by $\chi_{\textrm{II}}.$ Now we want to derive a system for $z=\int_0^{t} \chi_{\textrm{II}} (t-\tau) f(x(\tau)) d\tau$. By taking the time derivative, we get,
\begin{align*}
\dot{z}& = \chi(0) f(x(t)) + \int_0^{t} \dot{\chi}_{\textrm{II}} (t-\tau) f(x(\tau)) d\tau = A_0 f(x(t)) + \int_0^t (B_0 \chi_{\textrm{II}}+\chi_1)(t-\tau) f(x(\tau)) d\tau\\
&=A_0 f(x(t))  + B_0 z + z_1,
\end{align*}
where $z_1=\int_0^{t} \chi_1(t-\tau) f(x(\tau))d\tau$. Similarly, by taking derivative of $z_1$, we find,
\begin{align*}
\dot{z_1} &= \chi_1(0) f(x(t)) + \int_0^{t} \dot{\chi_1} (t-\tau)f(x(\tau))d\tau 
= A_1 f(x(t)) + \int_0^{t} B_1\chi_{\textrm{II}}(t-\tau) f(x(\tau)) d\tau\\
&= A_1 f(x(t))  + B_1 z.
\end{align*}

\subsection{The initial covariance matrix for the second order system (\ref{eq: n=2})}
We will derive the initial covariance matrix for system (\ref{eq: n=2}). The linear form of the system can be written as,
\begin{align}\label{secsys}
\dot{\left(\begin{array}{c}z_1 \\z\end{array}\right)} = 
D
\left(\begin{array}{c}z_1 \\z\end{array}\right)
+\left[\begin{array}{c}A_1 f +\xi_1\\A_0 f + \xi\end{array}\right],
\end{align}
and we assume the initial covariance matrix is:
\begin{align}
\left(\begin{array}{c}z_1(0) \\z(0)\end{array}\right) \left(\begin{array}{ccc}z_1(0) &z(0)\end{array}\right)
=Q,
\end{align}
where,
\begin{align}
D=\left[\begin{array}{ccc} 0 & B_1 \\ I & B_0\end{array}\right],\quad
Q=\left[\begin{array}{ccc} Q_1 & Q_{12} \\ Q_{12}^T & Q_2\end{array}\right].
\end{align}
We consider the case when $\xi_1$ and $\xi_2$ are uncorrelated. In order for the Lyapunov equation to hold, we need $DQ$ to be asymmetric, that is,
$$
-Q_2^T B_1^T=Q_1 + B_0 Q_{12}.
$$
Using linear differential equation (\ref{secsys}) and solving for the analytic solution directly, we have,
\begin{align*}
\left(\begin{array}{c}z_1(t) \\z(t)\end{array}\right)=
e^{Dt} \left(\begin{array}{c}z_1(0) \\z(0)\end{array}\right)
+\int_0^t e^{D(t-s)} \left(\begin{array}{c}\xi_1(s) \\\xi(s)\end{array}\right) ds
+\int_0^t e^{D(t-s)}\left(\begin{array}{c}A_1 \\A_0\end{array}\right) f(x(s)) ds.
\end{align*}
Then we  substitute it  back to $\dot{x}=z$, and  the equation is in the form of
\begin{align*}
\dot{x}=(0\quad I) e^{Dt} \left(\begin{array}{c}z_1(0) \\z(0)\end{array}\right)
+(0\quad I)\int_0^t e^{D(t-s)} \left(\begin{array}{c}\xi_1(s) \\\xi(s)\end{array}\right) ds
+\int_0^t (0\quad I)e^{D(t-s)}\left(\begin{array}{c}A_1 \\A_0\end{array}\right) f(x(s)) ds.
\end{align*}
In particular, the colored noise is given by,
\begin{align*}
w(t) = (0\quad I) e^{Dt} \left(\begin{array}{c}z_1(0) \\z(0)\end{array}\right)
+(0\quad I)\int_0^t e^{D(t-s)} \left(\begin{array}{c}\xi_1(s) \\\xi(s)\end{array}\right) ds,
\end{align*}
and thus the kernel function under this approximation can be expressed in terms of a matrix exponential,
\begin{align}\label{eq: chi2}
\chi_{\textrm{II}}(t) = (0\quad I)e^{Dt}\left(\begin{array}{c}A_1 \\A_0\end{array}\right).
\end{align}
Collecting terms, we have
\begin{align*}
\begin{split}
\langle w(t)w^T(t') \rangle & = (0 \quad I) e^{Dt} Q e^{D^Tt'} \left(\begin{array}{c}0 \\I\end{array}\right)
+\int_0^{t'} (0 \quad I) e^{D(t-s')} \Sigma e^{D^T (t'-s')}\left(\begin{array}{c}0 \\I\end{array}\right) ds'\\
&= (0 \quad I) e^{Dt} Q e^{D^Tt'} \left(\begin{array}{c}0 \\I\end{array}\right) 
+ (0 \quad I) e^{D(t-s')} Q e^{D^T(t'-s')} \left(\begin{array}{c}0 \\I\end{array}\right)\Big|_{0}^{t'}\\
&=(0 \quad I)e^{D(t-t')}Q\left(\begin{array}{c}0 \\I\end{array}\right).
\end{split}
\end{align*}
To match the FDT, we need,
\begin{align*}
\left[\begin{array}{ccc} Q_1 & Q_{12} \\ Q_{12}^T & Q_2\end{array}\right]\left(\begin{array}{c}0 \\I\end{array}\right)
=\left(\begin{array}{c}Q_{12} \\Q_2\end{array}\right) = \left(\begin{array}{c}A_1 \\A_0\end{array}\right),
\end{align*}
this leads us to the conditions,
$$
Q_{12}=A_1, \quad Q_2=A_0, \quad Q_1= -A_0^T B_1^T - B_0 A_1^T.
$$

With these choices, the FDT is exactly satisfied.

\subsection{The coefficients for the approximation with $n=2.$}  
The moments of second order rational approximation are
\begin{align}
\begin{split}
&[I-\lambda B_0-\lambda^2 B_1]^{-1} [\lambda A_0 +\lambda^2 A_1]\\
=& A_0 \lambda + ( B_0 A_0  + A_1)\lambda^2 +(B_0A_1 +B_0^2 A_0 + B_1 A_0)\lambda^3 +O(\lambda^4).
\end{split}
\end{align}

We first check the moment expansion for the kernel function, 
$$
\chi(\lambda)=[I+\lambda \gamma +\lambda \Theta(\lambda)]^{-1} \lambda.
$$
The Taylor expansion for the previous few terms is written as:
$$
\chi(\lambda)=M_0\lambda +M_1 \lambda^2 + M_2 \lambda^3 +O(\lambda ^4).
$$
Direct calculations yield,
\begin{align*}
&M_0=\frac{1-\lambda^2\Theta'(\lambda)}{[I+\lambda \gamma+\lambda \Theta(\lambda)]^2}(0)=I,\\
&M_1=-\frac{1}{2}\frac{2\gamma+2\Theta(\lambda) + 4 \Theta'(\lambda) \lambda + (I+\gamma)\Theta''(\lambda) \lambda^2
+ (\Theta''(\lambda)\Theta(\lambda)-2\Theta'(\lambda)^2)\lambda^3}
{[I+\lambda \gamma+\lambda \Theta(\lambda)]^3}(0)=-\gamma-\Theta(0),\\
&M_2=(\Theta(0)+\gamma)^2-\Theta'(0).
\end{align*}

Combining with the  expansion  for the rational function, we found,
\begin{align*}
&A_0=M_0,\quad B_0A_0+A_1 = M_1,\\
&B_0M_1+B_1A_0=M_2,\quad B_1=-A_1(\gamma +M_{\infty}).
\end{align*}
From this linear system, one easily finds that,
\begin{equation}
 B_0= M_2 (\gamma +M_{\infty})(I + \gamma +M_{\infty})^{-1}.
\end{equation}
The rest of the coefficients can be determined directly from these equations. 

\bigskip
For the particular example, we take, 
$$
\Theta(\lambda)=\frac{1}{2}\left(\sqrt{\frac{1}{\lambda^2}+4\omega_0^2}-\frac{1}{\lambda}\right),
$$
with expansion at $\lambda=0$  given by 
$$
K\lambda-\frac{K^2}{m}  \lambda^3 + O(\lambda^5).
$$
Therefore, we obtain
$$
\Theta(0)=0,\quad \Theta'(0)=K.
$$


\subsection{A short derivation of the Adelman and Doll model}
We start with the equations for the solid atoms, $m\ddot{x}_j = K(x_{j+1}-2x_j+x_{j-1}),$ {\blue with $x_j(0)=\dot{x}_j(0)=0$} for $j\ge 0.$
Taking the Laplace transform, we arrive at,
\begin{equation}
  X_{j+1} - \big(2 + m/K s^2\big) X_j + X_{j-1}=0, j\ge 0.
\end{equation}
The general solution can be written as $X_j = C \xi^j $, in which $\xi= 1 + \frac{m}{2K} s^2 - \frac{ms}{2K} \sqrt{s^2 + 4\frac{K}{m}}.$ The other root has modulus larger than one and it has to be ruled out since it leads to unbounded solutions. Thus we have $X_1=\xi X_0,$ and in the time domain, it is given by $x_1(t)=\int_0^t \beta(t-\tau) x_0(\tau) d\tau.$ To completely eliminate $x_1$, we consider the force on $x_0:$ 
$ K(x_1 - x_0).$ To proceed, we let $\theta(t)= K \int_t^{+\infty} \beta(\tau) d\tau.$ As a result, the Laplace transform is related as follows, $\Theta(s) =K (\xi(0) - \xi(s))/s= \frac{m}2 \sqrt{s^2 + 4\frac{K}{m}}-\frac{m}{2} s.$ In terms of the new kernel function, the force on $x_0$ is reduced to,
\begin{equation}
 K(x_1 - x_0) =\theta(t) x_0(0) - \int_0^t \theta(t-\tau) \dot{x}_0(\tau) d\tau.
\end{equation}

In the case when the initial conditions for $x_j$, $j>0$, is nonzero, an external force $f^\text{ex}(t)$ can be derived using the linear superposition principle. When the initial state is in thermal equilibrium, this force, together with $\theta(t) x_0(0)$, gives rise to the stationary Gaussian noise, i.e., $R(t)=f^\text{ex}(t) + \theta(t) x_0(0).$ This part of the calculation is less straightforward, and the readers are referred to \cite{LiE07} for a derivation using the Mori-Zwanzig formalism.

\bibliographystyle{plain}
\bibliography{GLERef,BD,other,GLERef1}

 \end{document}